\documentclass[twoside,11pt]{amsart}%{article}%
\usepackage{amsmath,amssymb,amsthm}

\begin{document}

\setlength{\textwidth}{145mm} \setlength{\textheight}{203mm}
\setlength{\parindent}{0mm} \setlength{\parskip}{2pt plus 2pt}

\frenchspacing

%\pagestyle{myheadings}

%\markboth{Dimitar Mekerov}{Some conditions the curvature tensor of
%the Bismut connection ... }

% THEOREM Environments ---------------------------------------------------
\numberwithin{equation}{section}
\newtheorem{thm}{Theorem}[section]
\newtheorem{lem}[thm]{Lemma}
\newtheorem{prop}[thm]{Proposition}
\newtheorem{cor}[thm]{Corollary}
\newtheorem{probl}[thm]{Problem}

\newtheorem{defn}{Definition}[section]
\newtheorem{rem}{Remark}[section]
\newtheorem{exa}{Example}%[section]

\newcommand{\be}[1]{\begin{equation}\label{#1}}
\newcommand{\ee}{\end{equation}}

% MATH -------------------------------------------------------------------

\newcommand{\X}{\mathfrak{X}}
\newcommand{\B}{\mathcal{B}}
\newcommand{\s}{\mathfrak{S}}
\newcommand{\g}{\mathfrak{g}}
\newcommand{\W}{\mathcal{W}}
\newcommand{\Lgr}{\mathrm{L}}
\newcommand{\dd}{\mathrm{d}}
\newcommand{\n}{\nabla}
\newcommand{\pd}{\partial}
\newcommand{\ddx}{\frac{\pd}{\pd x^i}}
\newcommand{\ddy}{\frac{\pd}{\pd y^i}}
\newcommand{\ddu}{\frac{\pd}{\pd u^i}}
\newcommand{\ddv}{\frac{\pd}{\pd v^i}}

\newcommand{\diag}{\mathrm{diag}}
\newcommand{\End}{\mathrm{End}}
\newcommand{\im}{\mathrm{Im}}
\newcommand{\id}{\mathrm{id}}

\newcommand{\ie}{i.e.}
\newfont{\w}{msbm9 scaled\magstep1}
\def\R{\mbox{\w R}}
\newcommand{\norm}[1]{\left\Vert#1\right\Vert ^2}
\newcommand{\nnorm}[1]{\left\Vert#1\right\Vert ^{*2}}
\newcommand{\nN}{\norm{N}}
\newcommand{\nJ}{\norm{\nabla J}}
\newcommand{\nnJ}{\nnorm{\nabla J}}
\newcommand{\tr}{{\rm tr}}

\newcommand{\thmref}[1]{Theorem~\ref{#1}}
\newcommand{\propref}[1]{Proposition~\ref{#1}}
\newcommand{\secref}[1]{\S\ref{#1}}
\newcommand{\lemref}[1]{Lemma~\ref{#1}}
\newcommand{\dfnref}[1]{Definition~\ref{#1}}

\frenchspacing

%%% ----------------------------------------------------------------------

\title[Connection with parallel totally skew-symmetric torsion ...]
{Connection with parallel totally skew-symmetric torsion %\\
on almost complex manifolds %\\
with Norden metric}

\author{Dimitar Mekerov}

%%% ----------------------------------------------------------------------
\maketitle
%%% ----------------------------------------------------------------------

{\small
%\begin{abstract}
\textbf{Abstract} \\
In the present work\footnote{Partially supported by project
RS09-FMI-003 of the Scientific Research Fund, Paisii Hilendarski
University of Plovdiv, Bulgaria} we consider an almost complex
manifold with Norden metric (i.~e. a metric with respect to which
the almost complex structure is an anti-isometry). On such a
manifold we study a linear connection preserving the almost
complex structure and the metric and having a totally
skew-symmetric torsion tensor. We consider the case when the
manifold admits a connection with parallel totally skew-symmetric
torsion and the case when such connection has a K\"ahler curvature
tensor. We get necessary and sufficient conditions for an
isotropic K\"ahler manifold with Norden metric. 
}
%\end{abstract}

%%%%%%%%%%%%%%%%%%%%%%%%%%%%%%%%%%%%%%%%%%%%%%%%%%%%%%%%%%%%%%%%%%%%%%%%%%%0

\section{Introduction}

There is a strong interest in the metric connections with totally
skew-symmet\-ric torsion tensor (3-form). These connections arise
in a natural way in theoretical and mathematical physics. For
example, such a connection is of particular interest in string
theory \cite{Stro}. In mathematics this connection was used by
Bismut to prove the local index theorem for non-K\"ahler Hermitian
manifolds \cite{Bis}. A connection with totally skew-symmetric
torsion tensor is called a KT connection by physicists,  and among
mathematicians this connection is known as a Bismut connection.

In the present work we continue the investigations from
\cite{Mek2} for a connection $\n'$ with totally skew-symmetric
torsion on non-K\"ahler quasi-K\"ahler manifolds with Norden
metric. There are proved some necessary and sufficient conditions
the curvature tensor of $\n'$ to be K\"ahlerian. In the case when
this tensor is K\"ahlerian, some relations between its scalar
curvature and the scalar curvatures of other curvature-like
tensors are obtained. Moreover, conditions for isotropic K\"ahler
manifolds with Norden metric are get.

Now we consider the case when $\n'$ has a parallel torsion. We
obtain a relation between the scalar curvatures for $\n'$ and the
Levi-Civita connection $\n$. We establish that the manifold is
isotropic K\"ahlerian with Norden metric iff these curvatures are
equal. We obtain a necessary and sufficient condition for $\n'$
with parallel torsion be with K\"ahler curvature tensor. Moreover
we show that if $\n'$ has a parallel torsion and a K\"ahler
curvature tensor, then the manifold is isotropic K\"ahlerian.

\section{Preliminaries}

Let $(M,J,g)$ be a $2n$-dimensional \emph{almost complex manifold
with Norden metric}, i.~e.
\begin{equation*}
J^2x=-x, \qquad g(Jx,Jy)=-g(x,y),
\end{equation*}
for all differentiable vector fields $x$, $y$ on $M$. The
\emph{associated metric} $\tilde{g}$ of $g$ on $M$, given by
$\tilde{g}(x,y)=g(x,Jy)$, is a Norden metric, too. The signature
of both metrics is necessarily $(n,n)$.

Further, $x$, $y$, $z$, $w$ will stand for arbitrary
differentiable vector fields on $M$ (or vectors in the tangent
space  of $M$ at an arbitrary point $p\in M$).

The Levi-Civita connection of $g$ is denoted by $\nabla$. The
tensor filed $F$ of type $(0,3)$ on $M$ is defined by
\begin{equation}\label{1.1}
F(x,y,z)=g\bigl( \left( \nabla_x J \right)y,z\bigr).
\end{equation}
It has the following properties:
\begin{equation*}
F(x,y,z)=F(x,z,y)=F(x,Jy,Jz),\quad F(x,Jy,z)=-F(x,y,Jz).
\end{equation*}

In \cite{GaBo}, the considered manifolds are classified into eight
classes with respect to $F$: $\W_0$, $\W_1$, $\W_2$, $\W_3$,
$\W_1\oplus\W_2$, $\W_1\oplus\W_3$, $\W_2\oplus\W_3$,
$\W_1\oplus\W_2\oplus\W_3$. The class $\W_0$ of the \emph{K\"ahler
manifolds with Norden metric} is contained in each of the other
seven classes. It is determined by the condition $F(x,y,z)=0$,
which is equivalent to $\n J=0$. The class
$\W_1\oplus\W_2\oplus\W_3$ is the class of all almost complex
manifolds with Norden metric.

The condition                                                                        %%%
\begin{equation}\label{1.2}
\mathop{\s} \limits_{x,y,z} F(x,y,z)=0,
\end{equation}
where $\mathop{\s} \limits_{x,y,z}$ is the cyclic sum over
$x,y,z$, characterizes the class $\W_3$ of the
\emph{quasi-K\"ahler manifolds with Norden metric}. This is the
only class among the basic classes $\W_1$, $\W_2$, $\W_3$ of
manifolds with non-integrable almost complex structure $J$.

Let $\{e_i\}$ ($i=1,2,\dots,2n$) be an arbitrary basis of the
tangent space of $M$ at a point $p\in M$. The components of the
inverse matrix of $g$, with respect to this basis, are denoted by
$g^{ij}$.

Following \cite{GRMa}, the \emph{square norm} $\nJ$ of $\nabla J$
is defined in \cite{MeMa} by
\begin{equation}\label{1.3}
    \nJ=g^{ij}g^{ks}
    g\bigl(\left(\nabla_{e_i} J\right)e_k,\left(\nabla_{e_j}
    J\right)e_s\bigr),
\end{equation}
where it is proven that
\begin{equation}\label{1.4}
    \nJ=-2g^{ij}g^{ks}
    g\bigl(\left(\nabla_{e_i} J\right)e_k,\left(\nabla_{e_s}
    J\right)e_j\bigr).
\end{equation}

There, the manifold with $\nJ=0$ is called an
\emph{isotropic-K\"ahler manifold with Norden metric}. It is clear
that every K\"ahler manifold with Norden metric is
isotropic-K\"ahler, but the inverse implication is not always
true.

Let $R$ be the curvature tensor of $\nabla$, i.~e. $
R(x,y)z=\nabla_x \nabla_y z - \nabla_y \nabla_x z -
    \nabla_{[x,y]}z$. The corresponding $(0,4)$-tensor is
determined by $R(x,y,z,w)=g(R(x,y)z,w)$. The Ricci tensor $\rho$
and the scalar curvature $\tau$ with respect to $\nabla$ are
defined by
\[
    \rho(y,z)=g^{ij}R(e_i,y,z,e_j),\qquad \tau=g^{ij}\rho(e_i,e_j).
\]

A tensor $L$ of type (0,4) with the pro\-per\-ties%
\be{1.5}%
L(x,y,z,w)=-L(y,x,z,w)=-L(x,y,w,z), \ee %
\be{1.6} %
\mathop{\s} \limits_{x,y,z} L(x,y,z,w)=0 \quad \textit{(the first
Bianchi identity)}
\ee %
is called a \emph{curvature-like tensor}. Moreover, if the
curvature-like tensor $L$ has the property
\begin{equation}\label{1.7}
L(x,y,Jz,Jw)=-L(x,y,z,w),
\end{equation}
it is called a \emph{K\"ahler tensor} \cite{GaGrMi}.

%%%%%%%%%%%%%%%%%%%%%%%%%%%%

Let $\n'$ be a linear connection with a tensor $Q$ of the
transformation $\n \rightarrow\n'$ and a torsion tensor $T$, i.~e.
\begin{equation}\label{1.8}
\n'_x y=\n_x y+Q(x,y),\qquad T(x,y)=\n'_x y-\n'_y x-[x,y].
\end{equation}
The corresponding (0,3)-tensors are defined by
\begin{equation}\label{1.9}
    Q(x,y,z)=g(Q(x,y),z), \qquad T(x,y,z)=g(T(x,y),z).
\end{equation}
The symmetry of the Levi-Civita connection implies
\begin{equation}\label{1.10}
    T(x,y)=Q(x,y)-Q(y,x), \qquad
    T(x,y)=-T(y,x).
\end{equation}

A linear connection $\n'$ on an almost complex manifold with
Norden metric $(M,J,g)$ is called a \emph{natural connection} if
$\n' J=\n' g=0$. The last conditions are equivalent to $\n' g=\n'
\tilde{g}=0$. If $\n'$ is a linear connection with a tensor $Q$ of
the transformation $\n \rightarrow\n'$ on an almost complex
manifold with Norden metric, then it is  a natural connection iff
the following conditions are valid:
\begin{equation}\label{1.12}
    F(x,y,z)=Q(x,y,Jz)-Q(x,Jy,z),
\end{equation}
\begin{equation}\label{1.13}
    Q(x,y,z)=-Q(x,z,y).
\end{equation}

According to \cite{Hay}, we have
\begin{equation}\label{1.14}
    Q(x,y,z)=\frac{1}{2}\bigl\{
    T(x,y,z)-T(y,z,x)+T(z,x,y)\bigr\}.
\end{equation}

Let $\n'$ be the natural connection with a totally skew-symmetric
torsion tensor $T$ on a non-K\"ahler manifold with Norden metric
$(M,J,g)$. According to \cite{Mek2} we have
\begin{equation}\label{2.12}
    Q(x,y)=\frac{1}{4}\bigl\{\left(\n_{x} J\right)Jy-\left(\n_{Jx} J\right)y-2\left(\n_{y} J\right)Jx
    \bigr\}.
\end{equation}
Since $\n'$ has a totally skew-symmetric torsion tensor $T$, then
\begin{equation}\label{2.2}
    T(x,y,z)=-T(y,x,z)=-T(x,z,y)=-T(z,y,x).
\end{equation}
From \eqref{1.14} and \eqref{2.2} it is follows that for the
tensor $Q$ it is valid
\begin{equation}\label{2.3}
    Q(x,y,z)=\frac{1}{2}T(x,y,z).
\end{equation}

%%%%%%%%%%%%%%%%%%%%%%%%%%%%%%%%%%%%---------------------------------------

\section{Connection with parallel totally skew-sym\-met\-ric
torsion}

Let $\n'$ be the connection with totally skew-symmetric torsion
tensor $T$ on the quasi-K\"ahler manifold with Norden metric
$(M,J,g)$.

Now we consider the case  when $\n'$ has a parallel torsion, i.~e.
$\n'T=0$.

It is known that the curvature tensors $R'$ and $R$ of $\n'$ and
$\n$, respectively,  satisfy
\begin{equation}\label{3.1}
    \begin{split}
    &R'(x,y,z,w)=R(x,y,z,w)+ \left(\n_x Q\right)(y,z,w)-\left(\n_y
    Q\right)(x,z,w)
    \\[4pt]
    &\phantom{K(x,y,z,w)=R(x,y,z,w)}+Q\bigl(x,Q(y,z),w\bigr)-Q\bigl(y,Q(x,z),w\bigr).
\end{split}
\end{equation}

Equality \eqref{2.3} implies $\n'Q=0$ in the considered case. Then
from the formula for covariant derivation with respect to $\n'$ it
follows that
\begin{equation}\label{3.2}
    xQ(y,z,w)-Q(\n'_x y,z,w)-Q(y,\n'_x z,w)-Q(y,z,\n'_x w)=0.
\end{equation}
According to the first equality of \eqref{1.8} we have
\begin{equation}\label{3.3}
    \begin{array}{l}
Q(\n'_x y,z,w)=Q(\n_x y,z,w)+Q(Q(x, y),z,w),\\[4pt]
Q( y,\n'_x z,w)=Q(y,\n_x z,w)+Q(y, Q(x, z),w),\\[4pt]
Q( y,z,\n'_x w)=Q(y,z,\n_x w)+Q(y,  z,Q(x,w)).\\[4pt]
    \end{array}
\end{equation}
Combining \eqref{3.2}, \eqref{3.3}, the first equality of
\eqref{1.9} and having in mind the formula for covariant
derivation with respect to $\n$, we obtain
\begin{equation}\label{3.4}
    \begin{split}
    &\left(\n_x Q\right)(y,z,w)=Q(Q(x, y),z,w)    \\[4pt]
    &\phantom{\left(\n_x Q\right)(y,z,w)=}
    -g\bigl(Q(x,z),Q(y,w)\bigr)-g\bigl(Q(y,z),Q(x,w)\bigr).
\end{split}
\end{equation}
From \eqref{3.4} and the first equality of \eqref{1.10} we have
\begin{equation}\label{3.5}
    \begin{split}
    &\left(\n_x Q\right)(y,z,w)-\left(\n_y Q\right)(x,z,w)
    =Q(T(x, y),z,w)    \\[4pt]
    &\phantom{\left(\n_x Q\right)(y,z,w)=}
    -2g\bigl(Q(x,z),Q(y,w)\bigr)+2g\bigl(Q(y,z),Q(x,w)\bigr).
\end{split}
\end{equation}

Because of \eqref{3.5}, equality \eqref{3.1} can be rewritten as
\begin{equation}\label{3.6}
    \begin{split}
    &R'(x,y,z,w)=R(x,y,z,w)+ Q\left(T(x,y),z,w\right)
    \\[4pt]
    &\phantom{R'(x,y,z,w)=}
    -g\bigl(Q(x,z),Q(y,w)\bigr)+g\bigl(Q(y,z),Q(x,w)\bigr).
    \end{split}
\end{equation}

Since $Q(e_i,e_j)=-Q(e_j,e_i)$ it follows that
$g^{ij}Q(e_i,e_j)=0$. Then, from \eqref{3.6} after contraction by
$x=e_i$, $w=e_j$, we obtain the following equality for the Ricci
tensor $\rho'$ of $\n'$:
\begin{equation}\label{3.7}
    \begin{split}
    &\rho'(y,z)=\rho(y,z)+2g^{ij}g\bigl(Q(e_i,y),Q(z,e_j)\bigr)
    \\[4pt]
    &\phantom{\rho'(y,z)=\rho(y,z)}
    -g^{ij}g\bigl(Q(e_i,z),Q(y,e_j)\bigr).
    \end{split}
\end{equation}
Contracting by $y=e_k$, $z=e_s$ in \eqref{3.7}, we get
\begin{equation}\label{3.8}
    \tau'=\tau+g^{ij}g^{ks}g\bigl(Q(e_i,e_k),Q(e_s,e_j)\bigr),
\end{equation}
where $\tau'$ is the scalar curvature of $\n'$.

By virtue of \eqref{3.8}, \eqref{2.12}, \eqref{1.3} and
\eqref{1.4} we have
\begin{equation}\label{3.9}
    \tau'=\tau-\frac{1}{8}\nJ.
\end{equation}

Thus we arrive at the following
\begin{thm}\label{thm-3.1}
Let $\n'$ be the connection with parallel totally skew-symmetric
torsion on the quasi-K\"ahler manifold with Norden metric
$(M,J,g)$. Then for the Ricci tensor $\rho'$ and the scalar
curvature $\tau'$ of $\n'$ are valid \eqref{3.7} and \eqref{3.9},
respectively. \hfill$\Box$
\end{thm}

Equality \eqref{3.9} leads to the following
\begin{cor}\label{cor-3.2}
Let $\n'$ be the connection with parallel totally skew-symmetric
torsion on the quasi-K\"ahler manifold with Norden metric
$(M,J,g)$. Then the manifold $(M,J,g)$ is isotropic K\"ahlerian
iff $\n'$ and $\n$ have equal scalar curvatures. \hfill$\Box$
\end{cor}

%%%%%%%%%%%%%%%%%%%%%%%%%%%%%%%%-------------------------------------------

\section{Connection with parallel totally skew-sym\-met\-ric torsion and K\"ahler curvature tensor}

Let $\n'$ be a connection with parallel totally skew-symmetric
torsion on the quasi-K\"ahler manifold with Norden metric
$(M,J,g)$.

We will find conditions for the curvature tensor $R'$ of $\n'$ to
be K\"ahlerian.

From \eqref{2.3}, having in mind that $Q$ is a 3-form, we have
\begin{equation*}
    \begin{split}
    &Q\left(T(x,y),z,w\right)=Q\left(z,w,T(x,y)\right)=g\bigl(Q(z,w),T(x,y)\bigr)
    \\[4pt]
    &\phantom{Q\left(T(x,y),z,w\right)}
    =g\bigl(T(x,y),Q(z,w)\bigr)=2g\bigl(Q(x,y),Q(z,w)\bigr).
    \end{split}
\end{equation*}

Then \eqref{3.6} obtains the form
\begin{equation}\label{4.1}
    \begin{split}
    &R'(x,y,z,w)=R(x,y,z,w)+2g\bigl(Q(x,y),Q(z,w)\bigr)
    \\[4pt]
    &\phantom{R'(x,y,z,w)=}
    -g\bigl(Q(x,z),Q(y,w)\bigr)+g\bigl(Q(y,z),Q(x,w)\bigr).
    \end{split}
\end{equation}

From \eqref{4.1}, identities \eqref{1.5} and \eqref{1.7} for $R'$
follow immediately. Therefore $R'$ is a K\"ahler tensor iff the
first Bianchi identity \eqref{1.6} for $R'$ is satisfied. Since
this identity is valid for $R$, then \eqref{4.1} implies that $R'$
is K\"ahlerian iff
\[
\mathop{\s} \limits_{x,y,z}
\bigl\{2g\bigl(Q(x,y),Q(z,w)\bigr)-g\bigl(Q(x,z),Q(y,w)\bigr)+g\bigl(Q(y,z),Q(x,w)\bigr)\bigr\}=0.
\]

Thus, using that $Q$ is a skew-symmetric tensor, we arrive the
following
\begin{thm}\label{thm-4.1}
Let $\n'$ be the connection with parallel totally skew-symmetric
torsion on the quasi-K\"ahler manifold with Norden metric
$(M,J,g)$. Then the curvature tensor for $\n'$ is a K\"ahler
tensor iff
\begin{equation}\label{4.2}
    \mathop{\s} \limits_{x,y,z}
g\bigl(Q(x,y),Q(z,w)\bigr)=0.
\end{equation} \hfill$\Box$
\end{thm}

Because of the skew-symmetry of $Q$, \eqref{4.2} implies
\begin{equation*}
    g\bigl(Q(y,z),Q(x,w)\bigr)-g\bigl(Q(x,z),Q(y,w)\bigr)=-g\bigl(Q(x,y),Q(z,w)\bigr).
\end{equation*}
The last equality and \eqref{4.1} lead to the following
\begin{cor}\label{cor-4.2}
Let $\n'$ be the connection with parallel totally skew-symmet\-ric
torsion and K\"ahler curvature tensor on the quasi-K\"ahler
manifold with Norden metric $(M,J,g)$. Then
\begin{equation}\label{4.3}
    R'(x,y,z,w)=R(x,y,z,w)+g\bigl(Q(x,y),Q(z,w)\bigr).
\end{equation}
\hfill$\Box$
\end{cor}

If $R'$ is a K\"ahler tensor then $R'(x,y,Jz,Jw)=-R'(x,y,z,w)$,
and because of \eqref{4.3} we have
\begin{equation}\label{4.4}
    \begin{split}
    &R(x,y,Jz,Jw)+R(x,y,z,w)=-g\bigl(Q(x,y),Q(Jz,Jw)\bigr)\\[4pt]
   &\phantom{R(x,y,Jz,Jw)+R(x,y,z,w)=}
    -g\bigl(Q(x,y),Q(z,w)\bigr).
    \end{split}
\end{equation}

From \eqref{2.12} we get
\[
Q(x,Jy)=JQ(x,y)-\left(\n_x J\right)y.
\]
Then we have
\[
Q(Jx,Jy)=-Q(x,y)-\left(\n_{Jx} J\right)y-\left(\n_{y} J\right)Jx
\]
and consequently
\[
    \begin{split}
&g\bigl(Q(x,y),Q(Jz,Jw)\bigr)=-g\bigl(Q(x,y),Q(z,w)\bigr)\\[4pt]
   &\phantom{g\bigl(Q(x,y),Q(Jz,Jw)\bigr)=}
-g\bigl(Q(x,y),\left(\n_{Jz} J\right)w+\left(\n_{w}
J\right)Jz\bigr).
    \end{split}
\]

The last equality and \eqref{4.4} imply the following
\begin{cor}\label{cor-4.3}
Let $\n'$ be the connection with parallel totally skew-symmetric
torsion and K\"ahler curvature tensor on the quasi-K\"ahler
manifold with Norden metric $(M,J,g)$. Then
\begin{equation}\label{4.5}
    R(x,y,Jz,Jw)+R(x,y,z,w)=g\bigl(Q(x,y),\left(\n_{Jz} J\right)w+\left(\n_{w}
J\right)Jz\bigr).
\end{equation}
\hfill$\Box$
\end{cor}

Contracting by $x=e_i$, $w=e_j$ in \eqref{4.5}, we obtain
\[
    g^{ij}R(e_i,y,Jz,Je_j)+\rho(y,z)=g^{ij}g\bigl(Q(e_i,y),\left(\n_{Jz} J\right)e_j+\left(\n_{e_j}
J\right)Jz\bigr).
\]
Then, after a contraction by $y=e_k$, $z=e_s$, it follows
\begin{equation}\label{4.6}
        \tau^{**}+\tau=g^{ij}g^{ks}g\bigl(Q(e_i,e_k),\left(\n_{Je_s} J\right)e_j+\left(\n_{e_j}
J\right)Je_s\bigr),
\end{equation}
where $\tau^{**}=g^{ij}g^{ks}R(e_i,e_k,Je_s,Je_j)$.

From \eqref{2.12}, \eqref{1.3} and  \eqref{1.4} we have
\[
g^{ij}g^{ks}g\bigl(Q(e_i,e_k),\left(\n_{Je_s}
J\right)e_j+\left(\n_{e_j} J\right)Je_s\bigr)=-\frac{1}{8}\nJ.
\]
Then \eqref{4.6} can be rewritten as
\[
\tau^{**}+\tau=-\frac{1}{8}\nJ.
\]
On the other hand, according to  \cite{MeMa}, we have
\[
\tau^{**}+\tau=-\frac{1}{2}\nJ.
\]
Then $\nJ=0$ and therefore the following is valid.
\begin{thm}
Let $(M,J,g)$ be  a quasi-K\"ahler manifold with Norden metric
which admit a connection with parallel totally skew-symmetric
torsion and K\"ahler curvature tensor. Then $(M,J,g)$ is a
isotropic K\"ahler manifold with Norden metric. \hfill$\Box$
\end{thm}

%%%%%%%%%%%%%%%%%%%%%%%%%%%%%%%%%%%%%%%%%%%%%%%%%%%%%%%%%%%%%%%%%%%%%%%%%%%

%%%%%%%%%%%%%%%%%%%%%%%%%%%%%%%%%%%%%%%%%%%%%%%%%%%%%%%%%%%%%%%%%%%%%%%%%%%%

\bigskip

\textit{Dimitar Mekerov\\
Department of Geometry\\
Faculty of Mathematics and Informatics
\\
Paisii Hilendarski University of Plovdiv\\
236 Bulgaria Blvd.\\
4003 Plovdiv, Bulgaria
\\
e-mail: mircho@uni-plovdiv.bg}

\end{document}